\documentclass[a4paper]{amsart}

\usepackage{amsthm,amsfonts,amsmath, amssymb}
  
\usepackage{xcolor}

\usepackage{setspace} 
\usepackage{anysize}
\marginsize{2.5cm}{3cm}{2cm}{2cm}   

\usepackage{graphicx,float,enumerate,subfigure,tikz}
\usepackage[ruled,boxed,commentsnumbered,norelsize]{algorithm2e}
\usepackage{verbatim}
\usepackage{url}
\usepackage[capitalise]{cleveref}

 \usepackage{enumitem}
\setlist[enumerate,1]{label=(\roman*),font=\normalfont}
 
\newtheorem{theorem}{Theorem}
\newtheorem{corollary}[theorem]{Corollary}
\newtheorem{lemma}[theorem]{Lemma}
\newtheorem{proposition}[theorem]{Proposition}
\newtheorem{problem}[theorem]{Problem}

\newenvironment{claimproof}{\noindent\textit{Proof.}}{\hfill$\square$}

\newtheorem{conjecture}[theorem]{Conjecture}

\theoremstyle{definition}

\theoremstyle{remark}

\crefname{remark}{Remark}{Remarks}

\theoremstyle{definition}
\newtheorem{claim}{Claim}
\newtheorem*{claim*}{Claim}
\usepackage{etoolbox}
\AtEndEnvironment{proof}{\setcounter{claim}{0}}

\theoremstyle{remark}

\usepackage{etoolbox}
\AtEndEnvironment{proof}{\setcounter{case}{0}}
\AtEndEnvironment{claimproof}{\setcounter{case}{0}}

\crefname{example}{Example}{Examples}

\usepackage{mathtools}

\DeclareMathOperator{\lk}{link}

\newcommand{\card}[1]{\ensuremath{\# #1}}

\newcommand{\R}{\mathbb{R}}

\newcommand{\E}{\mathcal{E}}

\newcommand{\Nc}{\mathcal{N}}

\newcommand{\V}{\mathcal{V}}

\crefname{rmk}{Remark}{Remarks}
\crefname{claim}{Claim}{Claims}
\crefname{problem}{Problem}{Problems}
\date{\today}
\title{Edge connectivity of simplicial polytopes}
 
\usepackage{amsaddr}

\author{Guillermo Pineda-Villavicencio \& Julien Ugon}
\address{Federation University Australia}
\address{School of Information Technology, Deakin University}
\email{\texttt{work@guillermo.com.au}}
\email{\texttt{julien.ugon@deakin.edu.au}} 

\thanks{Julien Ugon's research was partially supported by ARC discovery project DP180100602.}

\begin{document}
\begin{abstract} A simplicial polytope is a polytope with all its facets being combinatorially equivalent to simplices. We deal with the edge connectivity of the graphs of simplicial polytopes. We first establish that, 
for any $d\ge 3$, for any $d\ge 3$, every minimum edge cut of cardinality at most $4d-7$  in such a graph is \textit{trivial}, namely it consists of  all the edges incident with some vertex. A consequence of this is that, for $d\ge 3$,  
the graph of a simplicial $d$-polytope with minimum degree $\delta$  is $\min\{\delta,4d-6\}$-edge-connected.  In the particular case of $d=3$, we have that every minimum edge cut in a plane triangulation is trivial; this  may be of interest to researchers in graph theory. 

Second, for every $d\ge 4$ we construct a simplicial $d$-polytope whose graph has a nontrivial minimum edge cut of cardinality $(d^{2}+d)/2$. This gives  a  simplicial 4-polytope with a nontrivial minimum edge cut that has ten edges. Thus, the aforementioned result is best possible for simplicial $4$-polytopes.   
\end{abstract}   
\maketitle 

\section{Introduction} 
A (convex) polytope is the convex hull of a finite set $X$ of points in $\R^{d}$; the \textit{convex hull} of $X$ is  the smallest convex set containing $X$.  The \textit{dimension} of a polytope in $\R^{d}$ is one less than the maximum number of affinely independent points in the polytope; a set of points $p_{1},\ldots, p_{k}$ in $\R^{d}$ is {\it affinely independent} if  the $k-1$ vectors $p_{1}-p_{k},\ldots, p_{k-1}-p_{k}$ are linearly independent.  A polytope of dimension $d$ is referred to as a \textit{$d$-polytope}.

A polytope  is structured around other polytopes, its faces. A {\it face} of a polytope $P$ in $\R^{d}$ is $P$ itself, or  the intersection of $P$ with a hyperplane in $\R^{d}$ that contains $P$ in one of its closed halfspaces.  A face of dimension 0, 1, and $d-1$ in a $d$-polytope is a \textit{vertex}, an {\it edge}, and a {\it facet}, respectively. The sets of vertices and edges of a polytope or a graph are denoted by $\V$ and $\E$, respectively. The \textit{graph} $G(P)$ of a polytope  $P$ is the graph with vertex set $\V(P)$ and edge set $\E(P)$.   
 
Unless otherwise stated, the graph theoretical notation and terminology follow from \cite{BM08} and the polytope theoretical notation and terminology from \cite{Zie95}. Moreover, when referring to graph-theoretical properties of a polytope such as  minimum degree, connectivity, and edge connectivity, we mean properties of its graph.

This paper studies the edge connectivity of a simplicial polytope, namely the edge connectivity of the graph of the polytope. A {\it simplicial} polytope is a polytope with all its facets being simplices.  By a simplex we mean any polytope that is combinatorially isomorphic  to a simplex.  Two polytopes $P$ and $P'$ are \textit{combinatorially isomorphic} if their face lattices are isomorphic.

For a graph $G$, we say that a set $Z\subseteq \V(G)\cup \E(G)$ \textit{separates} two distinct vertices $x,y$, if every path in $G$ from $x$  to $y$ contains an element of $Z$. A  set $Z$ \textit{separates} $G$ if it separates two vertices of $G$.  
A graph is \textit{$r$-edge-connected}, for $r\ge 0$, if no two vertices  are separated by fewer than $r$ edges. And  a graph with at least $r+1$ vertices is \textit{$r$-connected} if removing any $r-1$ vertices leaves a connected subgraph.  Balinski \cite{Bal61} showed the following, which implies that the graph of a $d$-polytope is $d$-edge-connected.
    
 \begin{theorem}[Balinski {\cite{Bal61}}] For $d\ge 2$, the graph of a $d$-polytope is $d$-connected.
 \label{thm:Balinski}
 \end{theorem}
 
 We have recently studied other notions of connectivity in graphs of polytopes. The paper    \cite{ThiPinUgo18v3}  studied the vertex connectivity of cubical polytopes and showed  a cubical $d$-polytope with minimum degree $\delta$ is $\min\{\delta,2d-2\}$-connected, which implies that  $\min\{\delta,2d-2\}$-edge-connected. In addition, the paper \cite{BuiPinUgo20} analysed the linkedness of cubical polytopes, a stronger notion of connectivity.

A separating set $D$ of edges is an \textit{edge cut} in $G$ if there exists a nonempty proper subset $X\subseteq \V(G)$ such that $D= \E(X,\V(G)\setminus X)$; here $\E(X,\V(G)\setminus X)$ denotes the set of all edges from a vertex $x\in X$ to a vertex $y\in \V(G)\setminus X$. Henceforth,  we write $\E(X,\overline {X})$ instead of $\E(X,\V(G)\setminus X)$.   Every minimal separating set of edges in a connected graph $G$ is an edge cut \cite[Sec.~2.5]{BM08}. In addition, a \textit{trivial} edge cut is one whose edges are the ones incident with a single vertex; otherwise the edge cut is \textit{nontrivial}.
 We prove the following.

\vspace{0.2cm}

\noindent {\bf Theorem} (Edge connectivity theorem). {\it For every $d\ge 3$, in a simplicial $d$-polytope every minimum edge cut of cardinality at most $4d-7$ is trivial}.

\vspace{0.2cm}

  Let $\delta(G)$ or $\delta(P)$ be the minimum degree of a vertex in a graph  $G$ or a polytope $P$. A corollary trivially ensues.

\vspace{0.2cm}

\noindent {\bf Corollary}. {\it For every $d\ge 3$, a simplicial $d$-polytope with minimum degree $\delta$ is $\min\{\delta,4d-6\}$-edge-connected.}

\vspace{0.2cm}  

The graph of a 3-polytope is 3-connected and planar graph by Steinitz's theorem \cite{Ste22}, and a planar embedding of the graph of a simplicial 3-polytope is a \textit{plane triangulation} $G$; therefore, abusing terminology slightly, we will use interchangeably the terms  plane triangulation and simplicial 3-polytope.  In this case, Euler's formula \cite{Eul58,Eul58a} implies that minimum degree of $G$ is at most five.  Hence, a minimum separating set of edges in $G$ would have at most five edges.  Thus, another corollary is the following.
 
\vspace{0.2cm}

\noindent {\bf Corollary.} {\it Every minimum edge cut in a plane triangulation  is trivial.} 

\vspace{0.2cm}
For the interest of  graph theorists, we provide a short, graph-theoretical proof of this corollary  in \cref{sec:plane-edge-connectivity}; surprisingly, this result seems to be new.

Finally, for every $d\ge 4$, we construct a simplicial $d$-polytope with a minimum edge cut  that  contains $(d^{2}/2+d/2)$ edges and does not consist of the edges incident with a vertex.  This shows that the aforementioned theorem is best possible for simplicial 4-polytopes. We suspect the following is true.

\begin{conjecture} For every $d\ge 3$, there is a function $f(d)$ quadratic in $d$ such that, in a simplicial $d$-polytope, every minimum edge cut of cardinality at most $f(d)$ is trivial. 
\end{conjecture}

The  simplicial $d$-polytopes we constructed contain  an \textit{empty $(d-1)$-simplex}, a set of $d$ vertices  that does not form a face of the polytope but every proper subset does, and their graphs contain  large complete graphs.  So perhaps our theorem holds for all minimum edge cuts in certain classes of simplicial polytopes.

\begin{problem}\label{prob:bounds}
For $d\ge 3$, is every minimum edge cut of a flag $d$-polytope or a balanced $d$-polytope trivial?  
\end{problem} 

A \textit{flag polytope} is a simplicial polytope whose proper faces are the complete graphs of the graph, and so it contains no empty simplices. And a \textit{balanced $d$-polytope} is a  simplicial $d$-polytope whose vertices can be coloured with $d$ colours such that adjacent vertices receive different colours, which implies that its graph does not contain complete graphs with more than $d$ vertices.      

\section{Plane triangulations} 
\label{sec:plane-edge-connectivity}

Our first theorem (\cref{thm:plane-edge-connectivity}) is about plane triangulations and   may be of independent interest  to graph theorists, and so we provide an independent, graph-theoretical proof.   

For two vertices $x,y$ of a graph $G$, we say that a path $L:=x_{1}\ldots x_{n}$ is  an \textit{$x-y$ path} in $G$ if $\V(L)\cap \{x\}=\{x_{1}\}$  and $\V(L)\cap \{y\}=\{x_{n}\}$. 
For an edge cut $D$, let  $\V(D)$ denote the vertices in the edges of $D$.  

If $G$ is a  3-connected plane graph,  for a vertex $x$  let $\mathcal F_{x}$ be the faces of $G$ that contain $x$. The  \textit{link} of $x$ in $G$, denoted $\lk(x)$, is the subgraph of $G$ induced by the vertices and edges in $\mathcal F_{x}$ that are disjoint from $x$. It is a standard result that the link of a vertex in a 3-connected plane graph is a cycle that  contains the neighbours of the vertex; see, for instance, \cite[Cor.~10.8]{BM08}.  Since a plane triangulation can be considered as a simplicial 3-polytope, we could have used the standard definition of a link in \cref{sec:simplicial-edge-connectivity} and \cref{prop:link-polytope} to justify these assertions, but we want this section  to be devoid of polytope theory.

\begin{lemma} Let $G$ be a $d$-connected graph, and let $D:=\E(X,\overline{X})$ be a nontrivial minimum edge cut of $G$ for some $X\subseteq \V(G)$. Then $\#(X\cap \V(D))\ge d$ and $\#(\overline{X}\cap \V(D))\ge d$. 
\label{lem:edge-cut-small}
\end{lemma}
\begin{proof}
Let $n_X:=\#(X\cap \V(D))$ and let $\delta:=\delta(G)$.  If $n_X\leq d-1$, then removing these $n_X$ vertices does not disconnect the graph. Thus $X=X\cap \V(D)$. For any $x\in X$,  $x$ is incident with at least $\delta-n_X+1$ edges in $D$ (as $\deg(x) \geq \delta$). Since $\card D\le \delta$,  we have $n_X(\delta-n_X+1)\leq \delta$, and so $\delta\leq n_X\leq d-1$, a contradiction to the fact that $\delta\ge d$. Hence $n_X\ge d$. The same reasoning yields that $\#(\overline X\cap \V(D))\ge d$.
\end{proof}

\begin{theorem}[Edge connectivity of plane triangulations] Every minimum edge cut in a plane triangulation is trivial.
\label{thm:plane-edge-connectivity}
\end{theorem}
\begin{proof} Let $G$ be a plane triangulation, and let $D:=\E(X,\overline X)$ be a minimum  edge cut for some $X\subseteq \V(G)$.  By Euler's formula, the minimum degree of a plane graph is at most five. Thus  $\card D\le 5$. Suppose, by way of contradiction, that $D$ is nontrivial. It is a standard result of graph theory that a plane triangulation is 3-connected; this also follows from Balinski's theorem~\eqref{thm:Balinski}.  Let $n_{X}:=\#(X\cap \V(D))$ and $n_{\overline X}:=\#(\overline{X}\cap \V(D))$. Then, \cref{lem:edge-cut-small} ensures that  $n_{X}\geq 3$ and $n_{\overline X}\geq 3$.  

There is a vertex $u\in {X}\cap \V(D)$ with a unique neighbour $v$ in $\overline{X}\cap \V(D)$; if every vertex in $ {X}\cap \V(D)$ was incident with at least two edges from $D$, then $\card D\ge 6$, a contradiction. Moreover,  since $D$ is nontrivial, the vertex $u$ has at least one neighbour in $X$, say $u_1$. It follows that every $u_{1}-v$ path in $G$ passes through an edge of $D$. Thus, the edge cut $D$ must separate $u_{1}$ from $v$ in the link of $u$. Let $D_{u}:=D\cap \E(\lk(u))$.  The set $D_u$ contains precisely the edges incident with $v$ in $\E(\lk(u))$; otherwise either $D_{u}$ consists of the edges incident with $u_{1}$ or there would be at least two nonadjacent edges in $D_u$, both cases implying the existence of at least two neighbours of $u$ in $\overline{X}$. Because $\lk(u)$ is a cycle (see also \cref{prop:link-polytope}), we conclude that $v$ is incident with at least three edges in $D$, including $uv$.

Again, since $\card D\le 5$ and $v$ is incident with at least three edges in $D$, there is a vertex $v'\in \overline{X}\cap \V(D)$ with a unique neighbour $u'$ in ${X}\cap \V(D)$.  Following the same line of reasoning as before, the vertex $u'$ is incident with $u'v'$ and exactly two edges in $D\cap \lk(v')$, one of these two edges is $u'v$. By the 3-connectivity of $G$, the vertices $u', v$  cannot disconnect the graph. Thus,  there is at least one vertex in $(X\cap \V(D))\setminus \{u,u'\}$ adjacent to a vertex in $\overline{X}\setminus \{v',v\}$. Since this means that $\#D\geq 6$, we arrive at a contradiction. The proof of the theorem is complete.
\end{proof}

As a corollary we get the following.
\begin{corollary}A plane triangulation with minimum degree $\delta$ is $\delta$-edge-connected. 
\label{cor:plane-edge-connectivity}
\end{corollary}

\section{Simplicial polytopes}
\label{sec:simplicial-edge-connectivity}


 The \textit{boundary complex} of a polytope $P$ is the set of faces of $P$ other than $P$ itself. And the \textit{link} of a vertex $x$ in $P$, denoted $\lk(x)$, is the set of faces of $P$ that do not contain $x$ but lie in a facet of $P$ that contains $x$.  We require a result from \cite{Zie95}, which we proved in \cite[Prop.~12]{BuiPinUgo18}.   
   
\begin{proposition}[{\cite[Ex.~8.6]{Zie95}}]\label{prop:link-polytope} Let $P$ be a $d$-polytope. Then the link of a vertex in $P$ is combinatorially isomorphic to the boundary complex of a $(d-1)$-polytope. In particular, for each $d\ge 3$, the graph of the link of a vertex is isomorphic to the graph of a $(d-1)$-polytope.
\end{proposition}

%
%

Let $G$ be a graph and let $X\subseteq \V(G)$, then $G[X]$ denotes the  subgraph of $G$ induced by the set $X$ and $G-X$ denotes the subgraph  $G[\V(G)\setminus X]$ of $G$. We also require Menger's theorem~\cite{Menger1927}. 
 
 \begin{theorem}[Menger, {\cite{Menger1927}}]\label{thm:Menger}  
Let $G$ be a  graph, and let $x$ and $y$ be two nonadjacent of its vertices. Then the minimum number of vertices separating $x$ from $y$ in $G$ equals the maximum number of pairwise internally disjoint $x-y$ paths in $G$. 
\end{theorem}  

An extension of the reasoning in the proof of \cref{thm:plane-edge-connectivity} yields the main result of the paper. 

A vertex adjacent to a vertex $x$ in a graph $G$ is a \textit{neighbour} of $x$. We denote by $\Nc_{G}(x)$ the set of neighbours of $x$ in $G$.  We extend this notation to neighbours that belong to a subgraph or a subset of vertices of $G$; for instance, if $Y\subseteq \V(G)$  then $\Nc_{Y}(x)$ denote the neighbours of $x$ in $Y$.

\begin{theorem} 
For each $d\ge 3$,  in a  simplicial $d$-polytope  every minimum edge cut of cardinality at most $4d-7$   is trivial. 
\label{thm:simplicial-edge-connectivity}
\end{theorem}
\begin{proof} 
Let $P$ be a simplicial $d$-polytope, let $G$ be its graph,  let $\delta$ be the minimum degree of $G$, and  let $D:=\E(X,\overline X)$ be a minimum  edge cut of $G$ for some $X\subseteq \V(G)$.  By our hypothesis,  $\card D\le 4d-7$. The theorem holds for $d=3$ by \cref{thm:plane-edge-connectivity}, and so assume that $d\ge 4$. Let $n_{X}:=\#(X\cap \V(D))$ and $n_{\overline X}:=\#(\overline{X}\cap \V(D))$.  Suppose, by way of contradiction, that $D$ is nontrivial. By Balinski's theorem~\eqref{thm:Balinski}, $G$ is $d$-connected, in which case \cref{lem:edge-cut-small} ensures that  $n_{X}\geq d$ and $n_{\overline X}\geq d$. We need two simple claims.   

\begin{claim}  There is a vertex $w\in X\cap \V(D)$ such that $1\leq\card  \Nc_{\overline{X}}(w)\leq d-2$. Similarly,  there is a vertex $z\in \overline X\cap \V(D)$ such that $1\leq\card  \Nc_{{X}}(z)\leq d-2$.
\label{thm:simplicial-edge-connectivity-1}
\end{claim}
\begin{claimproof}  If, for every vertex $x$ in $X\cap \V(D)$, we have that $\card \Nc_{\overline{X}}(x)\geq d-1$, then $\#D \ge d(d-1) > 4d-7$ for $d\ge 4$ (as $n_{X}\ge d$), a contradiction.  The other statement is proved analogously.  
\end{claimproof}
  
\begin{claim} 
\label{thm:simplicial-edge-connectivity-2}
For every vertex $w\in X\cap \V(D)$ such that $1\leq \card \Nc_{\overline{X}}(w)\leq d-2$, then every vertex  in $\Nc_{\overline{X}}(w)$  is incident with at least $d +1- \card \Nc_{\overline{X}}(w)$ edges from $D$.  Similarly,  for every vertex $z\in \overline X\cap \V(D)$ such that $1\leq \card \Nc_{{X}}(z)\leq d-2$, then every vertex in $\Nc_{{X}}(z)$  is incident with at least $d +1- \card \Nc_{{X}}(z)$ edges from $D$. 
\end{claim}
\begin{claimproof} Since $D$ is nontrivial, there exists a vertex $w_1\in \Nc_X(w)$. Because $\card \Nc_{G}(w)\ge d$ and $\card \Nc_{\overline{X}}(w)\le d-2$, the vertex $w$ has neighbours in both $X$ and $\overline{X}$, and so there must exist a separating set of edges in $D\cap \E(\lk(w))$.  
  
Consider a vertex $v\in \Nc_{\overline{X}}(w)$. As the graph of $\lk(w)$ is isomorphic to the graph of a simplicial $(d-1)$-polytope (\cref{prop:link-polytope}), the vertex $v$ has degree at least $d-1$ in $\lk(w)$. Of the neighbours of $v$ in $\lk(w)$, at most $\card \Nc_{\overline X}(w)-1$ of them are in $\lk(w)\cap \overline X$, since $\V(\lk(w))=\Nc_{X}(w)\cup \Nc_{\overline X}(w)$. It follows that, of the neighbours of $v$ in $\lk(w)$,  at least $d-1 - (\card \Nc_{\overline{X}}(w)-1)$ are in $X$; in other words, $v$ is incident with at least $d -\Nc_{\overline{X}}(w)$ edges in $D\cap \E(\lk(w))$. If we also count the edge  $wv\notin \E(\lk(kw))$,  then we get the desired number of edges incident with $v$.  The statement about the vertex $z\in \overline X\cap \V(D)$  is proved analogously.    \end{claimproof}

Let $k: = \min (\{\card \Nc_{\overline{X}}(x): x\in X\cap \V(D)\}\cup  \{\card \Nc_{{X}}(x): x\in \overline X\cap \V(D)\})$, and let $u\in  \V(D)$ be a vertex such that $\card \Nc_{\overline{X}}(u)=k$. Without loss of generality, assume that $u\in X\cap\V(D)$.  Then, by  \cref{thm:simplicial-edge-connectivity-1} we have that $\card  \Nc_{\overline{X}}(u)\leq d-2$, and by  \cref{thm:simplicial-edge-connectivity-2}  every vertex $v$ in $\Nc_{\overline{X}}(u)$  is incident with $d+1-k$ edges of $D$. Every vertex in $\Nc_{\overline{X}}(u)$ is in $\V(D)\cap \overline X$, and so there are at least $d-k$ vertices in $(\V(D)\cap \overline X)\setminus \lk(u)$ (as $n_{\overline X}\ge d$).  Because every vertex in $\V(D)$ is incident with at least $k$ edges from $D$, for $k\ge 2$ we get that   
\begin{equation*}
	\#D \geq k(d+1-k) + k(d-k) = 2k(d-k) + k\geq 4d-6,
\end{equation*}
a contradiction. Therefore $k=1$. Denote by $v$ the unique vertex in $\Nc_{\overline{X}}(u)$. By  \cref{thm:simplicial-edge-connectivity-2} we have that  $\#\Nc_{X}(v)\geq d$.
 
If every vertex in $\overline X\cap \V(D)$ other than $v$ was incident with at least three edges in $D$, then we would have $\#D\geq d + 3(d-1) = 4d-3$, and so  there exists a vertex $v'\in \overline X\cap \V(D)$ such that $\card \Nc_{X}(v')\le 2$. We consider two cases according to the cardinality of $\Nc_{X}(v')$. 

First suppose that $\card \Nc_{X}(v')= 1$. Then the unique vertex $u'\in \Nc_{X}(v')$ is incident with at least $d$ edges from $D$ by~\cref{thm:simplicial-edge-connectivity-2}. Since $G$ is $d$-connected by  Balinski's theorem~\eqref{thm:Balinski},  $G$  would remain $(d-2)$-connected after removing $u'$ and $v$. By Menger's theorem~\eqref{thm:Menger}, there are  $(d-2)$ pairwise  internally disjoint $u-v'$ paths in $G-\{u',v\}$, and so $D$ contains $d-2$ pairwise disjoint edges not containing $u'$ or $v$, say $u_1v_1,\ldots,u_{d-2}v_{d-2}$ with $u_{i}\in X$ and $v_{i}\in \overline X$.    It is not possible that each $u_i$ is incident with at least three edges in $D$, as otherwise \cref{thm:simplicial-edge-connectivity-2} would ensure that every vertex $v_i$  ($i=1,\ldots,d-2$)  would be incident with at least $d-2$ edges in $D$. Thus, counting the edges in $D$ incident with $v$, $v'$, and $v_{1},\ldots, v_{d-2}$ we get that $\#D\geq d+1+3(d-2) = 4d-5$. So there is a vertex $u_i\in X$ adjacent to at most two vertices in $\overline{X}$, one of which is $v_i$. By~\cref{thm:simplicial-edge-connectivity-2}, this implies that $\#\Nc_X(v_i)\geq d+1-2$. Of the edges in $D$ not incident with any vertex in $\{v,v',v_i\}$, $d-3$ are edges $u_jv_j$ for $j\neq i$, and another $d-3$ of them are edges incident with $u'$. Therefore, \[\#D\geq \#\Nc_{X}(v) + \#\Nc_{X}(v_i) + \#\Nc_{X}(v') + 2(d-3) \ge d+ d-1 + 1 + 2(d-3) = 4d-6,\] another contradiction. Thus, $\card \Nc_{X}(v')= 2$.
    
The proof of the case $\card \Nc_{X}(v')= 2$ is analogous to the proof of the case $\card \Nc_{X}(v')= 1$, but we provide the details for the sake of completeness.  Let $\{u',u''\}:=\Nc_{X}(v')$. Then  \cref{thm:simplicial-edge-connectivity-2} yields that both $u'$ and $u''$  are incident with $d-1$ edges from $D$.  Removing $u'$, $u''$, and $v$  would make $G$ $(d-3)$-edge-connected by  Balinski's theorem~\eqref{thm:Balinski}. So $D$ contains $d-3$ pairwise disjoint edges not containing $u'$, $u''$, or $v$, say $u_1v_1,\ldots,u_{d-3}v_{d-3}$ with $u_{i}\in X$ and $v_{i}\in \overline X$.    It is not possible that each $u_i$ is incident with at least three edges in $D$; otherwise \cref{thm:simplicial-edge-connectivity-2} would ensure that every vertex $v_i, i=1,\ldots,d-3$ would be incident with at least $d-2$ edges in $D$, which yields that $\#D\geq d+2+3(d-2) = 4d-4>4d-7$ by counting the edges of $D$ incident with $v$, $v'$, and $v_{1},\ldots, v_{d-3}$. So there is a vertex $u_i\in X$ adjacent to at most two vertices in $\overline{X}$, one of which is $v_i$. 
By~\cref{thm:simplicial-edge-connectivity-2}, this implies that $\#\Nc_X(v_i)\geq d-1$. Of the edges in $D$ not incident with any vertex in $\{v,v',v_i\}$, $d-4$ are edges $u_jv_j$ for $j\neq i$, and another $2d-2-6$ of them are edges incident with $u'$ or $u''$. Therefore, for $d\ge 5$ we have that \[\#D\geq \#\Nc_{X}(v) + \#\Nc_{X}(v_i) + \#\Nc_{X}(v') + d-4+2d-8 \ge d+ d-1 + 2 + 3d-12 = 5d-11>4d-7,\] another contradiction. In the case $d=4$ we get $\card D=5d-11$. This means that $\#\Nc_{X}(v)=4$, $\#\Nc_{X}(v_i)=3$, and   $\#\Nc_{X}(v')=2$. But, since $n_{\overline X}\ge d$, there should exist a vertex $v''$ in $\overline X\cap \V(D)$ other than $\{v,v_{i},v'\}$; this gives an edge that has not been accounted for.  This final contradiction completes the proof of the theorem.    
\end{proof}

A simple corollary of \cref{thm:simplicial-edge-connectivity} is the following.

\begin{corollary} For $d\ge 3$, a simplicial $d$-polytope with minimum degree $\delta$ is
$\min\{\delta,4d-6\}$-edge-connected. 
\end{corollary}

\subsection{A construction of nontrivial minimum edge cuts} We construct a simplicial $d$-polytope that shows that the bound of $4d-7$ in \cref{thm:simplicial-edge-connectivity} is best possible for $d=4$. 

\begin{figure}  
\begin{center}     \includegraphics{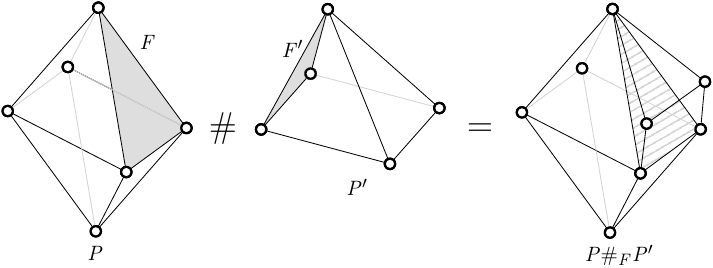}
\end{center}
\caption{Connected sum of two polytopes.}
\label{fig:connected-sum}  
\end{figure}

Let $P$ and $P'$  be two $d$-polytopes with a facet $F$ of $ P$  projectively isomorphic to a facet $F'$ of $P'$. The {\it connected sum} of  $P$ and $P'$ is obtained by ``gluing" $P$ and $P'$ along $F$ and $F'$.  Projective transformations on the polytopes $P$ and $P'$ may be required for the connected sum to be convex.   The connected sum of two polytopes is depicted in \cref{fig:connected-sum}. Our construction is based on performing connected sums along simplex facets. This is always possible because every polytope combinatorially isomorphic to a simplex is projectively isomorphic to the simplex \cite{McM76}.

We define stacked polytopes recursively. A $d$-simplex is \textit{stacked}, and each \textit{stacked} $d$-polytope with $n\ge d+2$ vertices is obtained as the connected sum of a stacked $d$-polytope with $n-1$ vertices and a $d$-simplex along a simplex facet.  

\begin{proposition} For each $d\ge 4$, there is a simplicial $d$-polytope with minimum degree at least  ${(d^2+d)}/{2}$ and a nontrivial minimum  edge  cut with ${(d^2+d)}/{2}$ edges. 
\label{prop:simplicial-nontrivial}
\end{proposition}
\begin{proof}
Consider a family of stacked $d$-polytopes $P_j$ with $j=d+1,\ldots, 2d$ vertices constructed as follows:
\begin{enumerate}
    \item $P_{d+1}$ is a $d$-simplex, with vertices labelled $x_1,\ldots, x_{d+1}$.
    \item For $j=d+1,\ldots, 2d-1$, the polytope $P_{j+1}$ is the connected sum of the $d$-simplex and the polytope $P_j$, along the simplex facet of $P_{j}$ containing the vertices $x_{j+1-d},x_{j+1-d+1},\ldots,x_{j}$. Label the last vertex added $x_{j+1}$. The vertices $x_{j+2-d},x_{j+2-d+1},\ldots,x_{j+1}$ form a facet of $P_{j+1}$, since they form a facet of the simplex we just added.
\end{enumerate}
The polytope $P_{2d}$ contains two disjoint simplex facets $F_0$ containing the vertices $x_1,\ldots,x_{d}$ and $F_1$ containing the vertices $x_{d+1},\ldots,x_{2d}$. Let  $D$ be the set of the edges in $P_{2d}$ from a vertex in $F_0$ to a vertex in  $F_1$. Then $D$ is a nontrivial edge cut of $P_{2d}$ with 
${(3d^2-d)}/{2} - d(d-1) = {(d^2+d)}/{2}$ edges.

Let $C$ be a cyclic $4$-polytope with degree at least ${(d^2+d)}/{2}$, and let $P$ be  the connected sum $C\#_{F_{0}} P_{2d}\#_{F_{1}}  C$. The polytope $P$ is the desired simplicial $d$-polytope. It has minimum degree at least ${(d^2+d)}/{2}$, and the set $D$ is a nontrivial edge cut of $P$. Say that $D=\E(X,\overline X)$ for some $X\subseteq \V(G)$. Then the edge cut $D$ is minimal in $P$, as both subgraphs $G[X]$ and $G[\overline X]$ are connected.

It remains to prove that $D$ is a minimum edge cut of $G(P)$. Let $C_{0}$ and $C_{1}$ be the two copies of $C$ in $P$. Then  $\E(P)=\E(C_{0})\cup \E(C_{1})\cup D$.  Consider a minimum edge cut $F$ in $G(P)$. It follows that $\card F\le \card D={(d^2+d)}/{2}$. We can assume that $F=F'\cup F''$ where $F'\subseteq D$ and $F''\cap D=\emptyset$. Because  $F''\subseteq \E(C_{0})\cup \E(C_{1})$ and both $C_{0}$ and $C_{1}$ are ${(d^2+d)}/{2}$-edge-connected, removing the edges of $F''$ does not disconnect the subgraph $G(C_{0})$ or $G(C_{1})$. In addition, if $F'$ is a proper subset of $D$, then removing $F'$ does not disconnect $G(P_{2d})$, as $D$ is minimal. It is now  plain that $F=D$, implying that $D$ is a minimum edge cut.

When $d=4$ this construction is best possible by \cref{thm:simplicial-edge-connectivity}. 
\end{proof}





\providecommand{\bysame}{\leavevmode\hbox to3em{\hrulefill}\thinspace}
\providecommand{\MR}{\relax\ifhmode\unskip\space\fi MR }
\providecommand{\MRhref}[2]{%
  \href{http://www.ams.org/mathscinet-getitem?mr=#1}{#2}
}
\providecommand{\href}[2]{#2}

\end{document}